\documentclass[10pt,twoside]{article}
\usepackage{epsf}
\usepackage{amsthm}
\usepackage{amssymb}
\usepackage{Latex-document}

\renewcommand{\thesection}{\arabic{section}}
\renewcommand{\thesubsection}{\arabic{section}.\arabic{subsection}}
\newcommand{\Section}[1]{\section{\hskip -1em.\hskip 1em#1}}
\newcommand{\Subsection}[1]{\subsection{\hskip -1em.\hskip 1em#1}}

\title{\bf Invariants of Legendrian Knots\vskip 6mm}

\author{Yu. V.  Chekanov\vspace*{-0.5cm}\thanks{Moscow Center for Continuous Mathematical Education,
B.~Vlasievsky per.~11, Moscow 119002, Russia.  E-mail: chekanov@mccme.ru}}

\date{\vspace{-8mm}}

\newtheorem{theorem}{Theorem}[section]

\newtheorem{corollary}[theorem]{Corollary}

\def \.{\mskip 1mu}
\def \?{\mskip -1mu}

\def \d{\partial}
\def \R{{\mathbb R}}
\def \N{{\mathbb N}}
\def \Z{{\mathbb Z}}

\def \im{{\rm im\.\.}}
\def \ker{{\rm ker\.\.}}
\def \ign#1{{}}
\def \Sw{{\rm Sw}}

\begin{document}

\maketitle

\thispagestyle{first} \setcounter{page}{385}

\begin{abstract}\vskip 3mm
We present two different constructions of
invariants for Legendrian knots in
the standard contact space $\R^3$.
These invariants are defined combinatorially,
in terms of certain planar projections, and
are useful in distinguishing Legendrian knots
that  have  the same classical invariants but
are not Legendrian isotopic.

\vskip 4.5mm


\end{abstract}

\vskip 12mm

\Section{Introduction} \label{intro}\setzero

\vskip-5mm \hspace{5mm}

\Subsection{Legendrian knots}

\vskip-5mm \hspace{5mm}

\begin{figure}[htb]
\centerline{\epsfbox{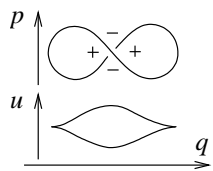}} \vskip -3mm \caption{Lagrangian projection and front projection} \label{1}
\end{figure}

A smooth knot $L$  in the standard contact space $(\R^3,\alpha)=(\{(q,p,u)\},\. du-p\.dq)$ is called Legendrian if
it is everywhere tangent to the $2$-plane distribution ${\rm ker }(\alpha)$ (or, in other words, if the
restriction of $\alpha$ to $L$ vanishes). Two Legendrian knots are Legendrian isotopic if the they can be
connected by a smooth path in the space of Legendrian knots (or, equivalently, if one can be sent to another by a
diffeomorphism $g$ of $\R^3$ such that $g^*\?\alpha=\varphi\alpha$, where $\varphi>0$). In order to visualize a
knot in $\R^3$, it is convenient to project it to a plane. In the Legendrian case, the character of the resulting
picture will depend on the choice of the projection. The useful two are: the Lagrangian projection
$\pi\colon\R^3\?\?\to\?\R^2,\,(q,p,u)\?\mapsto\? (q,p)$, and the front projection
$\sigma\colon\R^3\?\?\to\?\R^2,\,(p,q,u)\?\mapsto\? (q,u)$. In Figure~\ref{1}, two projections of the simplest
Legendrian knot (unknot) are shown.

We say that a Legendrian knot
$L\subset\R^3$ is $\pi$-generic if all self-intersections of the
immersed curve $\pi(L)$ are transverse double points.
We can represent a $\pi$-generic Legendrian knot $L$ by its
(Lagrangian) diagram:
the curve $\pi(L)\subset\R^2$, at every crossing of which
the overpassing branch (the one with the greater value of~$u$)
is marked.
Of course, not every abstract knot diagram in $\R^2$
is a diagram of a Legendrian knot, or is oriented
diffeomorphic to such
(it requires a bit of extra work to check whether
a given diagram corresponds to a Legendrian knot, cf.~\cite{C}).

Given a Legendrian knot $L\subset\R^3$, its
$\sigma$-projection, or front,
$\sigma(L)\subset\R^2$ is a singular curve
with nowhere vertical tangent vectors.
Its singularities, generically, are semi-cubic cusps and
transverse double points.
We say that $L$ is $\sigma$-generic if, moreover,
all self-intersections of $\sigma(L)$ have different $q$-coordinates.
Every closed planar curve with these types of singularities
and nowhere vertical tangent vectors is a front of a Legendrian knot.
Note that there is no need to explicitly indicate the
 type of a crossing of a front:
the overpassing branch (the one with the greater value of $p$)
is always the one with the greater slope.

\Subsection{Classical invariants}

\vskip-5mm \hspace{5mm}

The so-called classical
invariants of an oriented Legendrian knot $L$ are
defined as follows.
The first of them is, formally, just the smooth isotopy type of $L$.
The  Thurston--Bennequin number $\beta(L)$  of $L$
is the linking number (with respect to the orientation
defined by $\alpha\wedge d\alpha$)
between $L$ and $s(L)$,
where $s$ is a small shift along the $u$ direction.
The Maslov number $m(L)$ (which actually is an invariant
of Legendrian immersion) is twice the rotation number
of the projection of $L$ to the $(q,p)$ plane
(or, equivalently, the value of the Maslov
$1$-cohomology class on the fundamental class of $L$).
The change of orientation on $L$ changes the sign of $m(L)$
and preserves $\beta(L)$.
The Thurston-Bennequin number
of a $\pi$-generic Legendrian knot $L$ can be computed by counting
the crossings of its Lagrangian diagram $\pi(L)$ with signs:
$$
\beta(L)= \#\bigl(\raisebox{-1.5mm} {\makebox{\hskip0.3pt\epsfbox{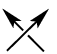}}}\hskip.3pt\bigr) \?
-\#\bigl(\raisebox{-1.5mm} {\makebox{\hskip0.3pt\epsfbox{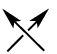}}}\hskip.3pt\bigr)
$$
(where the $q$ axis is horizontal and the $p$ axis is vertical).
In terms of the front projection,
the classical invariants can be computed as follows.
The Maslov number of a $\sigma$-generic oriented Legendrian knot $L$
is the number of  the right cusps of the front $\sigma(L)$, counted
with signs depending on the orientations:
$$
m(L)=\#\bigl(\raisebox{-1.5mm} {\makebox{\hskip0.5pt\epsfbox{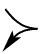}}}\hskip-.4pt\bigr)
\?-\#\bigl(\raisebox{-1.5mm} {\makebox{\hskip0.5pt\epsfbox{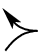}}}\hskip-.4pt\bigr).
$$
The Thurston-Bennequin
 number of $L$ is the number of crossings of $\sigma(L)$
counted with signs minus half the total number of cusps
(= the number of right cusps):
$$
\beta(L)= \#\bigl(\raisebox{-1.5mm} {\makebox{\hskip0.3pt\epsfbox{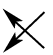}}}\hskip.2pt\bigr)
\?+\#\bigl(\raisebox{-1.5mm} {\makebox{\hskip0.3pt\epsfbox{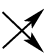}}}\hskip0.2pt\bigr)
\?-\#\bigl(\raisebox{-1.5mm} {\makebox{\hskip0.3pt\epsfbox{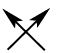}}}\hskip-.1pt\bigr)
\?-\#\bigl(\raisebox{-1.5mm} {\makebox{\hskip0.3pt\epsfbox{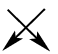}}}\hskip-.1pt\bigr)
\?-\#\bigl(\raisebox{-1.5mm} {\makebox{\hskip0.5pt\epsfbox{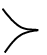}}}\hskip-.3pt\bigr).
$$
For the Legendrian knot shown in Figure~\ref{1}, we have $m=0$,
$\beta=-1$.

\Subsection{New invariants and classification results}

\vskip-5mm \hspace{5mm}

\begin{figure}[htb] \label{page-lag}
\centerline{\epsfbox{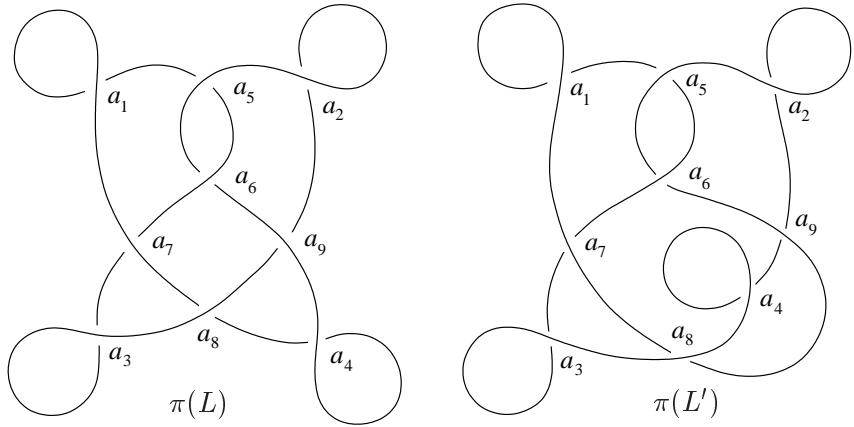}} \vskip -2mm \caption{Lagrangian projections of two Legendrian $5_2$ knots}
\label{2}
\end{figure}

\begin{figure}[htb]
\centerline{\epsfbox{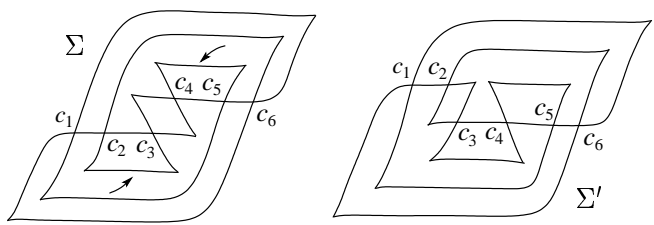}} \label{page-front} \vskip -2mm \caption{Fronts of two Legendrian $5_2$ knots}
\label{6}
\end{figure}

It is easy to show that every smooth knot admits a
Legendrian realization.
A natural question to ask is whether there exists a pair of
Legendrian knots which have the same classical invariants
but are not Legendrian isotopic.
The answer is positive:

\begin{theorem}{\hskip-2pt \rm\cite{C,CP}} \label{example}
There exist Legendrian knots $L, L'$
(see Figure~\ref{2} on p.~\pageref{page-lag},
Figure~\ref{6} on p.~\pageref{page-front})
that have the same classical invariants
(smooth knot type $5_2$,  $m=0$, $\beta=1$)
but are not Legendrian isotopic.
\end{theorem}
In the next two sections,
we present two combinatorial constructions of
Legendrian knot invariants.
The first one associates to the Lagrangian projection
of a Legendrian knot a differential graded algebra (DGA).
The second construction deals with decompositions of the
front projection into closed curves.
Each of the two  provides a proof of Theorem~\ref{example}.
The invariants do not change when the orientation of the knot reverses,
so essentially they are invariants of non-oriented Legendrian knots.
It should be mentioned that these  constructions
also produce, with minor modifications,
invariants of Legendrian links.

The number of Legendrian knots with given classical
invariants is known to be finite~\cite{CGH}.
Eliashberg and Fraser gave a classification of Legendrian
realization for smooth unknots~\cite{E,EF}.
It turned out that smooth unknots are
Legendrian simple in the sense that the Legendrian isotopy types
of their Legendrian realizations are determined by the classical
invariants.
Etnyre and Honda obtained a classification  of Legendrian
realization for torus knots and the figure eight knot~\cite{EH-cl}.
Again, these smooth knot types proved to be Legendrian simple.
The $5_2$ is the simplest knot type for which the
classification is not known.
By Theorem~\ref{example}, the type $5_2$ is not Legendrian simple.
Conjecturally, two Legendrian knots of smooth type $5_2$
with the same classical invariants are Legendrian isotopic
unless they form the  pair $L,L'$ from Theorem~\ref{example}.
Several interesting examples of knots with coinciding
classical invariants but not Legendrian isotopic were
constructed by Ng~\cite{N1,N2}.

\Section{DGA of a Legendrian knot} \label{pi} \setzero

\vskip-5mm \hspace{5mm}

\Subsection{Definitions}

\vskip-5mm \hspace{5mm}

In this section,
we associate with every $\pi$-generic Legendrian knot $L$
a DGA  $(A,\d)$ over $\Z/2\Z$
(\cite{C}; a similar construction was also given by Eliashberg).
This DGA is related to the symplectic field theory introduced
by Eliashberg, Givental, and Hofer in~\cite{EGH} (see~\cite{ENS}).

Let $\{a_1,\dots,a_n\}$
be  the  set of crossings of $Y=\pi(L)$.
Define $A$ to be the tensor algebra
(free associative unital algebra)
$T(a_1,\dots,a_n)$ with generators  $a_1,\dots,a_n$.
The grading on $A$ takes values
in the group $\Z/m(L)\Z$ and is defined as follows.
Given a crossing $a_j$, consider the points $z_+,z_-\in L$
such that $\pi(z_+)=\pi(z_-)=a_j$ and the $u$-coordinate of $z_+$
is greater than the $u$-coordinate of $z_-$.
These points
divide $L$ into two pieces, $\gamma_1$ and $\gamma_2$,
which we orient from $z_+$ to $z_-$.
We can assume, without loss of generality, that
the intersecting branches are orthogonal at $a$.
Then, for $\varepsilon\in\{1,2\}$, the
rotation number of  the curve $\pi(\gamma_\varepsilon)$
has the form $N_\varepsilon/2+1/4$,
where $N_\varepsilon\in\Z$.
Clearly, $N_1 - N_2$ is equal to $\pm m(L)$.
Hence $N_1$ and $N_2$
represent the same element of
the group  $\Gamma=\Z/m(L)\Z$,
which we define to be the degree of~$a_j$.

We are going to define the differential $\d$.
For every natural $k$, fix a (curved) convex $k$-gon
$\Pi_k\subset\R^2$ whose vertices
$x^k_0,\dots,x^k_{k-1}$ are numbered
counter-clockwise.
The form $dq\wedge dp$ defines an orientation on $\R^2$.
Denote by  $W_k(Y)$ the collection of
smooth orientation-preserving immersions
$f\colon\Pi_k\to\R^2$ such that $f(\d\Pi_k)\subset Y$.
Note that $f\in W_k(Y)$ implies
$f(x^k_i)\in\{a_1,\dots,a_n\}$.
Consider the set of nonparametrized immersions  $\widetilde W_k(Y)$,
which is the quotient of $W_k(Y)$
 by the action of the group
$\{g\in{\rm Diff_+} (\Pi_k)\.|\.g(x_i^k)=x_i^k\}$.
The diagram $Y$ divides a neighbourhood of each
of its crossings into four sectors.
We call positive two of them which are swept out by the underpassing curve
rotating counter-clockwise, and negative the other two
(the sectors are marked in Figure~\ref{1}).
For each vertex $x^k_i$ of the polygon $\Pi_k$,
a smooth immersion
$f\in \widetilde W_k(Y)$ maps its neighbourhood in $\Pi_k$
to either a positive or a negative sector;
we shall say that $x^k_i$ is, respectively,  a positive
or a negative vertex for $f$.
Define the set $W^+_k(Y)$ to consist of
immersions $f\in \widetilde W_k(Y)$ such that
the vertex $x^k_0$ is positive for $f$,
and all other vertices are negative.
Let $W^+_k(Y,a_j)=\{\.f\in W^+_k(Y)\,\. | \.\, f(x^k_0)= a_j\.\}$.
Denote $A_1=\{a_1,\dots,a_n\}\otimes\Z/2\Z\subset A$,
$A_k=(A_1)^{{\otimes}k}$.
Then $A=\oplus_{l=0}^\infty A_l $.
Let $\d=\sum_{k\ge 0} \d_k$,
where $\d_k(A_i)\in A_{i+k-1}$.
Define
$$
 \d_k(a_j)=\sum_{f\in W^+_{k+1}(Y,a_j)}
 f(x_1)\cdots f(x_k)
$$
(for $k=0$, we have
$\d_0(a_j)=\#( W^+_{1}(Y,a_j))$,
and extend $\d$ to $A$ by linearity and the Leibniz rule.
The following theorem says that $(A, \d)$ is indeed a DGA:

\begin{theorem} \label{main}
The differential $\d$ is well defined. We have $\deg(\d)\?=\?-1$
and $\d^2\?\?=\?0\?$.
\end{theorem}
Define the ($l$-th, where $l\?\in\? \Gamma$)
stabilization of a DGA $(T(a_1,\dots,a_n),\d)$
to be the DGA $(T(a_1,\dots,a_n, a_{n+1},a_{n+2}),\d)$, where
$\deg(a_{n+1})\?=\?l$, $\deg(a_{n+2})=l\?\?-\?\?1$, $\d(a_{n+1})=a_{n+2}$,
 and $\d$ acts on $a_1,\dots,a_n$ as before.
An automorphism of $T(a_1,\dots,a_n)$ is called elementary
if it sends $a_i$ to $a_i+v$, where $v$ does not involve $a_i$,
and fixes $a_j$ for $j\ne i$.
Two DGAs    $(T(a_1,\dots,a_n),\d)$, $(T(a_1,\dots,a_n),\d')$
are called tame isomorphic if one can be obtained from another by
a composition of elementary automorphisms;
they are called stable tame  isomorphic if they
become tame isomorphic after (iterated) stabilizations.
\begin{theorem} \label{move}
Let $(A,\d)$, $\,(A',\d')$ be the DGAs
of ($\pi$-generic) Legendrian knots $L,L'$.
If $L$ and $L'$ are Legendrian isotopic then
$(A,\d)$ and $\,(A',\d')$ are stable tame isomorphic.
In particular, the homology rings
$H(A,\d)=\ker(\d)/\im(\d)$ and $H(A',\d')=\ker(\d)/\im(\d)$
are isomorphic as graded rings.
\end{theorem}
The hard part in the proof of Theorem~\ref{main}
is to show that  $\d^2=0$.
The proof of this fact mimics, in a combinatorial way,
the classical gluing--compactness argument
of the Floer theory (cf.~\cite{F}).
The proof of Theorem~\ref{move} involves
a careful study of the behaviour of the DGA
associated with a Legendrian knot when its Lagrangian
diagram goes through elementary bifurcations
(Legendrian Reidemeister moves).

It turns out that one cannot replace the coefficient
ring $\Z/2\Z$ by $\Z\.$: in some sense, our homology theory is
not oriented.
However, the construction described above can be modified to
associate with a Legendrian knot $L$
a DGA graded by $\Z$ and having $\Z[s,s^{-\?1}]$
(where $\deg(s)=m(L)$) as a coefficient ring~\cite{ENS}.
After reducing the grading to $\Z/m(L)\Z$,
and applying the homomorphism $\Z[s,s^{-\?1}]\to\Z/2\Z $
sending both  $s$ and $1\?\in\?\Z$ to $1\in\?\Z/2\Z$,
this  $\Z[s,s^{-\?1}]\.$-DGA
becomes the $\Z/2\Z\.$-DGA of the knot $L$.

\Subsection{Poincar\'e polynomials}

\vskip-5mm \hspace{5mm}

Homology rings of DGAs can be hard to work with.
We are going to define an easily computable
invariant $I$, which is a finite subset of the group monoid
$\N_0[\Gamma]$, where $\N_0=\{0,1,\dots\}$, $\Gamma=\Z/m(L)\Z$.
Assume that $\d_0=0$. Then $\d_1^2=0$.
Since $\d(A_1)\subset A_1$, we can consider
the homology $H(A_1, \d_1)=\ker(\d_1|_{A_1})/\im(\d_1|_{A_1})$,
which is a vector space graded by the cyclic group $\Gamma$.
Define the Poincar\'e polynomial
$P_{(\?A,\d)}\?\?\in\N_0[\Gamma]$ by
$$
P_{(\?A,\d)}(t)=\sum_{\lambda\in\Gamma}
\dim\bigl(H_\lambda(A_1,\d_1)\bigr)\.t^\lambda,
$$
where $H_\lambda(A_1,\d_1)$ is the degree $\lambda$
homogeneous component of $H(A_1, \d_1)$.
Define the group ${\rm Aut}_0(A)$ to consist of graded
automorphisms of $A$ such that for each $i\in\{1,\dots,n\}$
we have $g(a_i)=a_i+c_i$, where $c_i\in A_0=\Z/2\Z$.
(of course, $c_i=0$ when $\deg(a_i)\ne 0$).
Consider the set $U_0(A,\d)$ consisting of automorphisms
$g\in {\rm Aut}_0(A)$ such that $(\d^g)_0\?=0$
(where $\d^g\?\?=g^{-1}\?\?\circ\d\circ g$).
Define
$$
I(A,\d) = \{ \.\.P_{(\?A,\d^{\.g})}\. \.\.|\.\.\.
g\in U_0\?(A,\d)\.\.\}.
$$
Since ${\rm Aut}_0(\?A)$ has at most $2^n$ elements,
this invariant is not hard to compute.
We can associate with every ($\pi$-generic)
Legendrian knot $L$ the set $I(L)=I(A_L,\d_L)$.
Note that  $P(-1)=\beta(L)$ for $P\in I(L)$.
One can show that $I$ is an invariant of stable tame DGA  isomorphism.
Hence Theorem~\ref{move} implies the following
\begin{corollary}\label{I}
If $L$ is Legendrian isotopic to $L'$ then $I(L)=I(L')$.
\end{corollary}
The set $I(L)$ can be empty (cf.~Section~\ref{inst})
but no examples are known where $I(L)$ contains more than
one element.
Also, for all known examples of pairs $L,L'$ of Legendrian knots with
coinciding classical invariants we have $P(1)=P'(1)$,
where  $P\?\in\? I(L)$,  $P'\?\in\? I(L')$.
Other, more complicated invariants of stable tame isomorfism
were developed and applied to distinguishing Legendrian knots
in~\cite{N2}.

\Subsection{Examples} \label{examples}

\vskip-5mm \hspace{5mm}

\noindent {\bf 1a.}\enspace Let  $(A,\d)=(T(a_1,\dots,a_9),\d)$ be the DGA of the Legendrian knot L given in
Figure~\ref{2}. We have  $m(L)=0$, $\beta(L)=1$,
 $\,\deg(a_i)=1$ for $i\le4$, $\,\deg(a_5)=2$, $\,\deg(a_6)=-2$,
 $\,\deg(a_i)=0$ for $i\ge7$,
$\,\d(a_1)= 1 +a_7 + a_7 a_6 a_5,$
$\,\d(a_2)= 1+ a_9 + a_5 a_6 a_9,$
$\,\d(a_3)= 1 +  a_8 a_7,$
$\,\d(a_4)= 1 +  a_8 a_9,$
$\,\d(a_i)=0$ for $i\ge5$.

\noindent
{\bf 1b.}\enspace
Let $(A',\d)=(T(a_1,\dots,a_9),\d)$
be the DGA
of the Legendrian knot $L'$ given in Figure~\ref{2}.
We have  $m(L')=0$, $\beta(L')=1$,
$\,\deg(a_i)=1$ for $i\le4$,
 $\,\deg(a_i)=0$ for $i\ge5$,
$\,\d(a_1)= 1 +a_7 + a_5+ a_7 a_6 a_5 + a_9 a_8 a_5,$
$\,\d(a_2)= 1+ a_9 + a_5 a_6 a_9,$
$\,\d(a_3)= 1 +  a_8 a_7, $
$\,\d(a_4)= 1 +  a_8 a_9, $
$\,\d(a_i)=0$ for $i\ge5$.

An explicit computation shows that
$I(L)=\{t^{-2}\?+t^1\?+t^2\}$,
$\.I(L')=\{2t^0\?+t^1\}$ and hence
Theorem~\ref{example} follows from Corollary~\ref{I}.

\begin{figure}[htb]
\centerline{\epsfbox{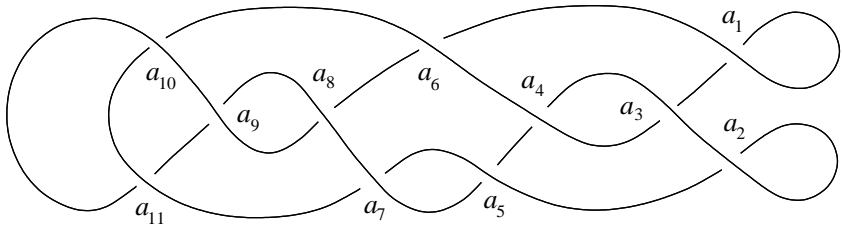}} \vskip -1.5mm \caption{Lagrangian projection of a Legendrian $6_2$ knot}
\label{3}
\end{figure}

\noindent
{\bf 2.}\hskip2pt\cite{N1} \enspace
Let $(A,\d)=(T(a_1,\dots,a_{11}),\d)$
be the DGA of the Legendrian knot $K$ given in Figure~\ref{3}.
We have  $m(K)=0$, $\beta(K)=-7$,
$\,\deg(a_i)=1$ for $i\in\{1,2,7,9,10\}$,
 $\,\deg(a_i)=0$ for $i\in\{3,4\}$,
$\,\deg(a_i)=-1$ for $i\in\{5,6,8,11\}$;
$\,\d(a_1)= 1 + a_{10} a_5 a_3,$
$\,\d(a_2)= 1+ a_3 + a_3 a_6 a_{10} + a_3 a_{11} a_7,$
$\,\d(a_4)= a_5 +  a_{11} +  a_{11} a_7 a_5,$
$\,\d(a_6)= a_{11} a_8,$
$\,\d(a_7)= a_8 a_{10} ,$
$\,\d(a_9)= 1 + a_{10}a_{11},$
$\,\d(a_i)=0$ for $i\in\{3,5,8,10,11\}$.
Denote by $\widehat K$ the `Legendrian mirror' of $K$ --- the image
of $K$ under the map
$(q,p,u)\?\mapsto\?(-\?q,p,-\?u)$.
The Legendrian knots $K, \widehat K$ have the same
classical invariants.
However, they are not Legendrian isotopic, and it is possible
to distinguish them by means of their DGAs.
There exist homology classes
$\xi_+, \xi_-$ in the graded homology ring $H(A,\d)$ such that
$\deg(\xi_+)=1,$ $\deg(\xi_-)=-1$, and $\xi_+\.\. \xi_-\?=1$
(choose $\xi_+\?\?=\?[a_{10}],$ $\xi_-\?\?=\?[a_{11}]$).
It follows from the definitions that
the DGA for $\widehat K$ is obtained from  $(A,\d)$
by applying the anti-automorphism reversing the order of
generators in all monomials.
Thus, if $K$ and $\widehat K$ are Legendrian isotopic then
the graded homology ring $H(A,\d)$ is anti-isomorphic to
itself, and   there exist
$\xi'_+, \xi'_-\in H(A,\d)$ such that $\deg(\xi_+')=1$,
$\deg(\xi'_-)=-1$,
$\.\xi'_-\.\.\xi'_+\? =1$.
But one can check that such classes do not exist
(see~\cite{N1,N2} for details) and hence $K$ and $\widehat K$
are not Legendrian isotopic.
Note that `first order invariants' such as  Poincar\'e polynomials
are useless in distinguishing Legendrian mirror knots.

\Section{Admissible decompositions of fronts}\label{sigma}

\setzero\vskip-5mm \hspace{5mm}

\Subsection{Definitions}

\vskip-5mm \hspace{5mm}

In this section, we present the
invariants of Legendrian knots constructed in~\cite{CP}.
These invariants are defined in terms of the front projection.

\begin{figure}[htb]
\centerline{\epsfbox{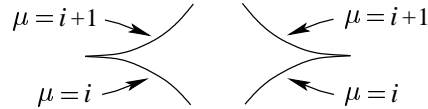}} \vskip -2mm \caption{Jumps of the Maslov index near cusps} \label{4}
\end{figure}

Given a $\sigma$-generic oriented Legendrian knot $L$,
denote by $C(L)$ the set of its points corresponding
to cusps of $\sigma(L)$.
The Maslov index $\mu\colon L\setminus C(L)\to \Gamma=\Z/m(L)\Z$
is a locally constant function, uniquely defined up to an additive
constant by the following rule:
the value of $\mu$ jumps at points of $C(L)$
by $\pm 1$ as shown in Figure~\ref{4}.
We call a crossing of $\Sigma=\sigma(L)$ Maslov if $\mu$ takes
the same value on both its branches.

\begin{figure}[htb]
\centerline{\epsfbox{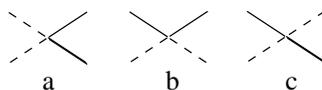}} \vskip -2.5mm \caption{Local decompositions} \label{5}
\end{figure}

Assume that $\Sigma=\sigma(L)$ is a union of closed
curves $X_1, \dots, X_n$ that have finitely many
self-intersections and meet each other at
finitely many points.
Then we call the unordered collection  $\{X_1, \dots, X_n\}$
a decomposition of $\Sigma$.
A decomposition $\{X_1, \dots, X_n\}$ is called admissible
if it satisfies certain conditions, which we are going to define.
The first two are as follows:
\begin{itemize}
\item[(1)]
Each curve $X_i$ bounds a topologically embedded disk: $X_i=\partial B_i$.
\item[(2)]
For each $i\in\{1,\dots,n\}$, $q\in\R$,
the set  $B_i(q)=\{\,u\in\R\mid(q,u)\in B_i\.\}$
is either a segment, or consists
of a single point $u$ such that $(q,u)$ is a cusp of $\Sigma$,
or is empty.
\end{itemize}
Conditions (1) and (2) imply that each curve $X_i$
has exactly two cusps (and hence the number of curves
is half the number of cusps).
Each $X_i$ is divided by cusps into two pieces,
on which the coordinate $q$ is a monotone function.
Near a crossing $x\in X_i\cap X_j$, the
decomposition of $\Sigma$
may look in one of the  three ways represented
in Figure~\ref{5}.
Conditions (1) and (2), in particular, rule out
the decomposition shown in Figure~\ref{5}a.
We call the crossing point $x$ switching
if $X_i$ and $X_j$ are not smooth near $x$
(Figure~\ref{5}b), and non-switching otherwise
(Figure~\ref{5}c).
\begin{itemize}
\item[(3)]
 If $(q_0,u)\in X_i\cap X_j$ ($i\ne j$)
is switching then for each $q\ne q_0$ sufficiently
close to $q_0$ the set $B_i(q)\cap B_j(q)$ either
coincides with  $B_i(q)$ or $B_j(q)$, or is empty.
\item[(4)] Every switching crossing is Maslov.
\end{itemize}
We call a decomposition admissible if it satisfies Conditions (1)-(3),
and graded admissible if it also satisfies Condition (4).
Denote by ${\rm{Adm}}(\Sigma)$
(resp.~${\rm{Adm}}_+(\Sigma)$) the set of admissible
(resp.~graded admissible)
decompositions of $\Sigma$.
Given $D\in{\rm{Adm}}(\Sigma)$,
denote by ${\Sw(D)}$ the set of its  switching points.
 Define $\theta(D)=\#(D))-\#(\Sw(D))$.
\begin{theorem}\label{adm}
If $\sigma$-generic Legendrian knots  $L,L'\subset\R^3$ are
Legendrian isotopic then there exists a one-to-one
mapping $g\colon {\rm{Adm}}(\sigma(L))\to {\rm{Adm}}(\sigma(L'))$
such that $g({\rm{Adm}}_+(\sigma(L)))= {\rm{Adm}}_+(\sigma(L'))$ and
$\theta (g(D))=\theta(D)$ for each  $D\in {\rm{Adm}}(\sigma(L))$.
In particular, the numbers $\# ({\rm{Adm}}(\sigma(L)))$ and
 $\# ({\rm{Adm}}_+(\sigma(L)))$
are invariants of Legendrian isotopy.
\end{theorem}

\Subsection{Remarks}

\vskip-5mm \hspace{5mm}

{\bf 1.}\enspace  Decompositions of fronts were first
considered by Eliashberg in~\cite{E87}
(only  Conditions (1) and (2) were involved).

\medskip\noindent
{\bf 2.}\enspace
No examples are known where the total number of admissible
decompositions  $\# ({\rm{Adm}}(\Sigma))$ is
different for two Legendrian knots with coinciding classical
invariants.

\medskip\noindent
{\bf 3.}\enspace
The proof of Theorem~\ref{adm} goes as follows:
we connect $L$ with $L'$ by a generic path
in the space of Legendrian knots and
define a canonical way to extend
admissible decompositions through
the points where the front is not $\sigma$-generic.
The mapping $g$ depends on the choice of the path:
a loop in the space of Legendrian knots
lifts to an automorphism
of  ${\rm{Adm}}(\sigma(L))$ which can be non-trivial
even when the loop is contractible.
The meaning of this phenomenon is not clear.

\medskip\noindent
{\bf 4.}\enspace It would be interesting to understand the
relation between admissible decompositions and DGAs of
Legendrian knots.
The first result in this direction is that
if ${\rm{Adm}}_+(\sigma(L))$ is nonempty then
the set $I(L)$ defined in the previous section is
also nonempty~\cite{Fu}.

\Subsection{Examples}

\vskip-5mm \hspace{5mm}

Note that every admissible decomposition $D$ of
a front $\Sigma$
is uniquely defined by its set of switching points.
Indeed, denote by $X(\Sigma)$ the set of crossings of $\Sigma$,
then each subset $E\subset X(\Sigma)$
defines a decomposition $D(E)$ of $\Sigma$ which
near $x\in X(\Sigma)$ has the form shown in Figure~\ref{5}b
if $x\in E$, and the form shown in Figure~\ref{5}c otherwise.
Clearly, if $E=\Sw(D)$ then $D=D(E)$.

The Legendrian knots represented by the
fronts $\Sigma,\Sigma'$ in
Figure~\ref{6} are respectively Legendrian isotopic
to the Legendrian knots $L,L'$ defined in Figure~\ref{2}.
We are going to show that
$\# ({\rm{Adm}}_+(\Sigma))=1$, $\# ({\rm{Adm}}_+(\Sigma'))=2$,
and hence Theorem~\ref{example}
is a consequence  of Theorem~\ref{adm}.
Assume that $D\in {\rm{Adm}}(\Sigma)$.
Consider the curve $X_1\in D$ containing
the piece of $\Sigma$ indicated by the lower arrow.
Being applied to $X_1$, Conditions (1) and (2)
imply that $c_2,c_3\in \Sw(D)$.
Similarly, looking at the curve $X_2\in D$
containing the piece of $\Sigma$
indicated by the upper arrow,
we conclude that  $c_4,c_5\in \Sw(D)$.
If one of the crossings $c_1, c_6$ is switching,
so is the other. Then either
$\Sw(D)=\{c_2,c_3,c_4,c_5\}$ or
 $\Sw(D)=\{c_1,c_2,c_3,c_4,c_5,c_6\}$.
It is not hard  to check that both decompositions
are admissible but only the first one is graded.
Thus  $\# ({\rm{Adm}}_+(\Sigma))=1$.
Arguing similarly, one can find that $\# ({\rm{Adm}}(\Sigma'))=2$,
where the admissible decompositions $D_1,D_2$
are defined by
 $\Sw(D_1)=\{c_2,c_3,c_4,c_5\}$,
 $\Sw(D_2)=\{c_1,c_2,c_3,c_4,c_5,c_6\}$, and are both graded.

\Section{Instability of invariants}\label{inst}

\setzero\vskip-5mm \hspace{5mm}

\begin{figure}[htb]
\centerline{\epsfbox{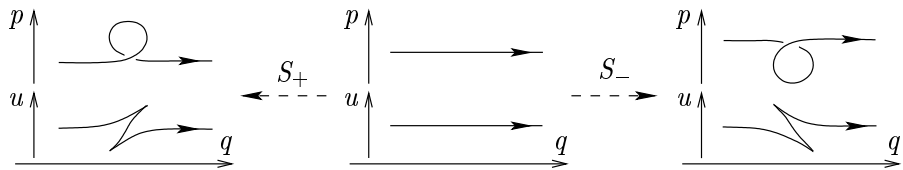}} \vskip -2mm \caption{Stabilizations} \label{7}
\end{figure}

There are two stabilizing operations,
$S_-$ and $S_+$, on Legendrian isotopy classes
of oriented Legendrian knots, defined as follows.
Given an oriented Legendrian knot $L$,
we perform one of the operations shown in Figure~\ref{7}
in a small neighbourhood of a point on $L$.
One can check that, up to Legendrian isotopy,
the resulting Legendrian knot $S_\pm(L)$ does not depend
on the choices involved, and the operations
$S_-,S_+$ commute.
An important observation is that
two Legendrian knots $L,L'$ have the same classical invariants
if and only if they are stable Legendrian isotopic in the sense that
there exist $n_-,n_+\in \N_0$ such that  $S_-^{n_-}(S_+^{n_+}(L))$
is Legendrian isotopic to $S_-^{n_-}(S_+^{n_+}(L'))$ \cite{FT}.

Thus the invariants constructed in Sections~\ref{pi} and~\ref{sigma}
cannot be stable.
In fact, they fail already after the first stabilization.
The homology ring $H$ of the DGA corresponding to $S_\pm(L)$ vanish,
and the  set $I(S_\pm(L))$ is empty.
This can be easily derived from the fact that the DGA of $S_\pm(L)$
can be obtained from the DGA of $L$ by adding a new generator
$a$ such that $\d(a)=1$.
The front of  $S_\pm(L)$ has no admissible decompositions because
Conditions (1) and (2) cannot hold
for the curve $X_i$ containing the newly created cusps.

Studying Legendrian realizations of non-prime knots,
Etnyre and Honda constructed, for each $m$, examples of Legendrian knots
that have the same classical invariants but are not Legendrian isotopic
even after $m$ stabilizations~\cite{EH-sum}.
Their proof uses the classification of Legendrian torus knots given
in~\cite{EH-cl}.
It is an open problem to find invariants distinguishing
those knots, or any pair of stabilized knots with the same
classical invariants.

\label{lastpage}

\end{document}